\documentclass[12pt]{amsart}
\usepackage{amsfonts, amsmath, amssymb, latexsym}
\usepackage[all]{xy}
	\title[Homological Epimorphisms of DGAs]{Homological Epimorphisms of Differential Graded Algebras}
    \author{David Pauksztello}
    \date{11th September 2006}    

\newtheorem{dfn}{Definition}[section]

\newtheorem{thrm}[dfn]{Theorem}

\newtheorem{exa}[dfn]{Example}

\newtheorem{lem}[dfn]{Lemma}

\newtheorem{prp}[dfn]{Proposition}

\newtheorem{corr}[dfn]{Corollary}

\newtheorem{rmk}[dfn]{Remark}

\newtheorem{rmks}[dfn]{Remarks}

\newtheorem{set}[dfn]{Setup}

\newtheorem{quest}[dfn]{Question}

\newtheorem{answ}[dfn]{Answer}

\newenvironment{prf}{\noindent\textbf{Proof: }}{$\Box$\ }

\newcommand{\comment}[1]{}

\newcommand{\D}{\mathcal{D}}

\newcommand{\K}{\mathcal{K}}

\newcommand{\F}{\mathcal{F}}

\newcommand{\Fop}{\mathcal{F}^{\rm{op}}}

\newcommand{\Rhom}[3]{\textnormal{RHom}_{#1}(#2,#3)}

\newcommand{\Ltensor}[1]{\stackrel{L}{\otimes}_{#1}}

\newcommand{\rightiso}{\stackrel{\sim}{\longrightarrow}}

\newcommand{\qtn}[1]{\begin{quest}\emph{#1}}

\newcommand{\ans}[1]{\begin{answ}\emph{#1}}

\newcommand{\eqtn}{\end{quest}}

\newcommand{\eans}{\end{answ}}

\newcommand{\df}[1]{\begin{dfn}\emph{#1}}

\newcommand{\edf}{\end{dfn}}

\newcommand{\thm}{\begin{thrm}}

\newcommand{\ethm}{\end{thrm}}

\newcommand{\lemma}{\begin{lem}}

\newcommand{\elemma}{\end{lem}}

\newcommand{\prop}{\begin{prp}}

\newcommand{\eprop}{\end{prp}}

\newcommand{\rem}[1]{\begin{rmk}\emph{#1}}

\newcommand{\erem}{\end{rmk}}

\newcommand{\rems}[1]{\begin{rmks}\emph{#1}}

\newcommand{\erems}{\end{rmks}}

\newcommand{\pf}{\begin{prf}}

\newcommand{\epf}{\end{prf}}

\newcommand{\cor}{\begin{corr}}

\newcommand{\ecor}{\end{corr}}

\newcommand{\setup}[1]{\begin{set}\emph{#1}}

\newcommand{\esetup}{\end{set}}

\newcommand{\eg}[1]{\begin{exa}\emph{#1}}

\newcommand{\eeg}{\end{exa}}

\entrymodifiers={+!!<0pt,\fontdimen22\textfont2>}

\begin{document}
    
\address{Department of Pure Mathematics, University of Leeds,
Leeds. LS2 9JT, United Kingdom}
\email{davidp@maths.leeds.ac.uk}
\urladdr{http://www.maths.leeds.ac.uk/\~{ }davidp}

\begin{abstract}
Let $R$ and $S$ be differential graded algebras. In this paper we give a characterisation of when a differential graded $R$-$S$-bimodule $M$ induces a full embedding of derived categories
$${}_{R}M_{S}\Ltensor{S}-:\D(S)\rightarrow \D(R).$$
In particular, this characterisation generalises the theory of Geigle and Lenzing's homological epimorphisms of rings, described in \cite{Geigle}. Furthermore, there is an application of the main result to Dwyer and Greenlees's Morita theory.
\end{abstract}

\maketitle

\setcounter{section}{-1}
	
\section{Introduction}

In \cite{Geigle}, Geigle and Lenzing characterise, using the classical derived functors $\textnormal{Ext}$ and $\textnormal{Tor}$, when a homomorphism of rings $\phi: R\rightarrow S$ induces a full embedding of bounded derived categories, $\D^{b}(S)\hookrightarrow \D^{b}(R)$. Geigle and Lenzing refer to such ring homomorphisms as homological epimorphisms of rings.

Differential graded algebras (DGAs) may be regarded as a generalisation of rings. It is, therefore, natural to ask whether the characterisation of Geigle and Lenzing also works for DGAs. It turns out that it does. Moreover, the characterisation for DGAs is a special case of a more general result regarding differential graded bimodules (DG bimodules). 

Given two DGAs $R$ and $S$ and a DG $R$-$S$-bimodule $M$ we can look at the functor $${}_{R}M_{S}\Ltensor{S}-:\D(S)\rightarrow \D(R)$$ and ask when this is a full embedding of derived categories.

The case of homological epimorphisms of DGAs then becomes the situation when $M=S$, with $S$ acquiring the left $R$-structure via a morphism of DGAs $\phi: R\rightarrow S$. In the work of Keller (\cite{Keller2}, Remarks 3.2), DG bimodules are regarded as generalised morphisms of DGAs. Thus, asking when ${}_{R}M_{S}\Ltensor{S}-$ is a full embedding of derived categories is analogous to asking when $M$ is a generalised homological epimorphism of DGAs. This more general setting makes the structural reasons behind the characterisation of homological epimorphisms of rings in \cite{Geigle} more transparent.

This paper is organised in the following way. In Section \ref{prelim} we give a brief exposition of derived functors, compact objects and canonical maps, which will be used in later sections. Section \ref{new} contains the main results of this paper, characterising when a DG $R$-$S$-bimodule $M$ induces a full embedding of derived categories $${}_{R}M_{S}\Ltensor{S}-:\D(S)\rightarrow \D(R).$$ In Section \ref{examples} we consider three examples of the main result. The first example relates to Dwyer and Greenlees's Morita theory of \cite{Greenlees}. The second example is the characterisation of homological epimorphisms of rings given by Geigle and Lenzing in \cite{Geigle}. The third example is a generalisation of Geigle and Lenzing's characterisation of homological epimorphisms of rings to DGAs.

It is worth emphasising again that the motivation of the results in Section \ref{new} is to give clear structural reasons why the characterisation of homological epimorphisms of rings by Geigle and Lenzing exists. 

\section{Preliminaries}\label{prelim}

We begin by reviewing the definitions of DGAs and DG modules, which can be found in \cite{Bernstein}.

\df{A \textit{differential graded algebra} $R$ over the commutative ground ring $\mathbb{K}$ is a graded algebra $R=\bigoplus_{i\in\mathbb{Z}}R_{i}$ over $\mathbb{K}$ together with a differential, that is, a $\mathbb{K}$-linear map $\partial^{R}:R\rightarrow R$ of degree -1 with $\partial^{2}=0$, satisfying the Leibnitz rule
$$\partial^{R}(rs)=\partial^{R}(r)s+(-1)^{|r|}r\partial^{R}(s)$$
where $r,s\in R$ and $r$ is a graded element of degree $|r|$.}
\edf

\df{A \textit{differential graded left $R$-module} (\textit{DG left $R$-module}) $M$ is a graded left module $M=\bigoplus_{i\in\mathbb{Z}}M_{i}$ over $R$ (viewed as a graded algebra) together with a differential, that is, a $\mathbb{K}$-linear map $\partial^{M}:M\rightarrow M$ of degree -1 with $\partial^{2}=0$, satisfying the Leibnitz rule
$$\partial^{M}(rm)=\partial^{R}(r)m+(-1)^{|r|}r\partial^{M}(m)$$
where $m\in M$ and $r\in R$ is a graded element of degree $|r|$. DG right $R$-modules are defined similarly. We shall denote the category of DG left $R$-modules by DG-Mod$(R)$.}
\edf

We denote by $R^{\rm{op}}$ the opposite DGA of $R$, that is, $R^{\rm{op}}$ consists of the same underlying set but the multiplication is given by $r\cdot s=(-1)^{|r||s|}sr$ for $r,s\in R$. DG right $R$-modules can be canonically identified with DG left $R^{\rm{op}}$-modules. Thus, we will denote the category of DG right $R$-modules by DG-Mod$(R^{\rm{op}})$. Henceforth, ``$M$ is a DG $R$-module'' will mean that $M$ is a DG left  $R$-module and ``$M$ is a DG $R^{\rm{op}}$-module'' will mean that $M$ is a DG right $R$-module. $M$ is said to be a DG $R$-$S$-bimodule if it is a DG $R$-module and a DG $S^{\rm{op}}$-module, with the $R$ and $S^{\rm{op}}$ structures compatible, that is, $r(ms)=(rm)s$ for all $r\in R$, $s\in S$ and $m\in M$. Subscripts are used to emphasise whether $M$ is a DG left or right module. 

A morphism $f:M\rightarrow N$ of DG $R$-modules is an $R$-linear map which is also a morphism of complexes of degree zero. Two morphisms $f,g:M\rightarrow N$ of DG $R$-modules are said to be \textit{homotopic} if there exists a homomorphism of graded modules $h:M\rightarrow N$ of degree $+1$ such that $$f-g=\partial^{N}h+h\partial^{M}.$$ Note that homotopy is an equivalence relation. The \textit{homotopy category} of $R$, written $\K(R)$, is the category whose objects are DG $R$-modules and whose morphisms are homotopy classes of morphisms of DG $R$-modules. The \textit{derived category} of $R$, $\D(R)$, is given by formally inverting the quasi-isomorphisms in $\K(R)$ (that is, the morphisms which induce isomorphisms in (co)homology). Recall that $\K(R)$ and $\D(R)$ are triangulated categories. Since $\D(R)$ is triangulated, it is equipped with an autoequivalence $\Sigma:\D(R)\rightarrow \D(R)$, called the \textit{suspension functor}, which is defined as follows. Let $M$ be a DG $R$-module, then:
$$(\Sigma M)^{n}=M^{n+1}, \ \textrm{and} \ \partial^{\Sigma M}=-\partial^{M}.$$
See \cite{Hartshorne} for details.

We shall now give brief details of compact objects, derived functors and canonical maps.

\subsection{Compact objects}\label{compact}

Let $R$ be a DGA. We take the following definition from \cite{Recollement}.

\df
{A DG $R$-module $M$ is \textit{finitely built from} ${}_{R}R$ in $\D(R)$ if $M$ can be obtained from ${}_{R}R$ using finitely many distinguished triangles, suspensions, direct summands and finite coproducts.}
\edf

\df
{A DG $R$-module $C$ is called \textit{compact} in $\D(R)$ if the functor $\textnormal{Hom}_{\D(R)}(C,-)$ commutes with set indexed coproducts in $\D(R)$.}
\edf

\rem
{It can be shown that a DG $R$-module $M$ is finitely built from ${}_{R}R$ in $\D(R)$ if and only if it is compact in $\D(R)$ (see \cite{Neeman}).}
\erem

\subsection{Derived functors}\label{derived}

Let $R$ be a DGA. For DG $R$-modules $M$ and $N$, we can construct a  complex $\textnormal{Hom}_{R}(M,N)$; similarly, for a DG $R^{\rm{op}}$-module $M$ and a DG $R$-module $N$ we can construct a complex $M\otimes_{R} N$. These are analogues of the classical homological functors Hom and tensor. There exist derived functors of these functors analogous to the classical derived functors, denoted $\Rhom{R}{M}{N}$ and $M\Ltensor{R}N$, see \cite{Hartshorne}. We make a brief note on how these are computed. 

We recall the following definitions from \cite{Bernstein}.

\df
{A DG $R$-module $P$ is said to be \textit{K-projective} if the functor $\textnormal{Hom}_{R}(P,-)$ preserves quasi-isomorphisms. A DG $R$-module $I$ is said to be \textit{K-injective} if the functor $\textnormal{Hom}_{R}(-,I)$ preserves quasi-isomorphisms. A DG $R$-module $F$ is said to be \textit{K-flat} if $F\otimes_{R}-$ preserves quasi-isomorphisms. Similarly for DG $R^{\rm{op}}$-modules.}
\edf

\df
{Let $M$ be a DG $R$-module. A \textit{K-projective resolution} of $M$ is a K-projective DG $R$-module $P$ together with a quasi-isomorphism $\pi : P\rightarrow M$. A \textit{K-flat resolution} of $M$ is a K-flat DG $R$-module $F$ together with a quasi-isomorphism $\theta : F\rightarrow M$. A \textit{K-injective resolution} of $M$ is a K-injective DG $R$-module $I$ together with a quasi-isomorphism $\iota: M\rightarrow I$.}
\edf 

Let $M$ and $N$ be DG $R$-modules. We compute $\Rhom{R}{M}{N}$ in one of the following two ways. First, we obtain a K-projective resolution $\pi: P\rightarrow M$ and note that $\Rhom{R}{M}{N}\cong \Rhom{R}{P}{N}\cong \textnormal{Hom}_{R}(P,N)$, or we obtain a K-injective resolution $\iota : N\rightarrow I$ and note that $\Rhom{R}{M}{N}\cong \Rhom{R}{M}{I}\cong \textnormal{Hom}_{R}(M,I)$. The derived functor $M\Ltensor{R}N$ is computed by substituting by a K-projective or a K-flat resolution of $M$ or $N$.

\subsection{Canonical maps}\label{canonical}

When considering the classical homological functors of Hom and tensor we have the notion of a ``canonical'' map. For example, if $R$ and $S$ are rings and $M$ is a left $S$-module then there is a natural map ${}_{S}M\rightarrow \textnormal{Hom}_{R}({}_{R}S_{S},{}_{R}S_{S}\otimes_{S}{}_{S}M)$ given by $m\mapsto (s\mapsto s\otimes m)$. A similar construction can be found for tensor. These natural maps of the classical homological functors induce natural maps of the classical derived functors Ext and Tor. However, in the setting of DGAs where the analogous ``hyper-homological'' derived functors occur as the right derived Hom functor, RHom, and the left derived tensor functor, $\Ltensor{}$, which are computed by substitution by K-projective or K-injective resolutions, we have to be more careful in constructing ``canonical'' maps.

However, the derived functors $\textnormal{RHom}$ and $\Ltensor{}$ are adjoints and we may use adjunction to find a suitable candidate for a ``canonical'' map. For example, consider the natural map ${}_{S}M\rightarrow \textnormal{Hom}_{R}({}_{R}S_{S},{}_{R}S_{S}\otimes_{S}M)$ given above. There is a natural isomorphism
$$\textnormal{Hom}_{R}({}_{R}S_{S}\otimes_{S}{}_{S}M,{}_{R}S_{S}\otimes_{S}{}_{S}M) \cong \textnormal{Hom}_{S}({}_{S}M,\textnormal{Hom}_{R}({}_{R}S_{S},{}_{R}S_{S}\otimes_{S}{}_{S}M))$$
given by adjunction. It is easy to see that the natural map $m\mapsto (s\mapsto s\otimes m)$ corresponds to the identity map on ${}_{R}S_{S}\otimes_{S}{}_{S}M$ under the adjunction.

We give an example of the corresponding situation with $R$ and $S$ as DGAs and $N$ as a DG $S$-module, which we shall use immediately in Section \ref{new}. Let $M$ be a DG $R$-$S$-bimodule. We may ask: what is the canonical map ${}_{S}N\rightarrow \Rhom{R}{{}_{R}M_{S}}{{}_{R}M_{S}\Ltensor{S}{}_{S}N}$ for a DG $S$-module $N$?

Here, as above, we have the following isomorphism given by adjunction:
\begin{align*}
\lefteqn{\textnormal{Hom}_{\D(R)}({}_{R}M_{S}\Ltensor{S}{}_{S}N,{}_{R}M_{S}\Ltensor{S}{}_{S}N)} & \\
\hspace{5ex} & \cong \textnormal{Hom}_{\D(S)}({}_{S}N,\Rhom{R}{{}_{R}M_{S}}{{}_{R}M_{S}\Ltensor{S}{}_{S}N}).
\end{align*}
Then we define the canonical map ${}_{S}N\rightarrow \Rhom{R}{{}_{R}M_{S}}{{}_{R}M_{S}\Ltensor{S}{}_{S}N}$ to be the image of the identity map  on ${}_{R}M_{S}\Ltensor{S}{}_{S}N$ under the adjunction isomorphism.

\section{Main results}\label{new}

Recall that for a functor $F:\mathcal{C}\rightarrow \D$ between two categories $\mathcal{C}$ and $\D$ the terms \textit{fully faithful} and \textit{full embedding} are synonymous and will be used interchangeably; the second usage arising because of the fact that $F$ being fully faithful means that, up to equivalence, $\mathcal{C}$ can be regarded as a full subcategory of $\D$.

\prop
\label{new1}
Let $R$ and $S$ be DGAs and suppose that $M$ is a DG $R$-$S$-bimodule. Then the following conditions are equivalent:

$(1)$ For all DG $S$-modules $N$ the canonical map $${}_{S}N\rightarrow \Rhom{R}{{}_{R}M_{S}}{{}_{R}M_{S}\Ltensor{S}{}_{S}N}$$ is an isomorphism.

$(2)$ For all DG $S$-modules $N$ and $N'$ the canonical map $$\Rhom{S}{{}_{S}N}{{}_{S}N'}\rightarrow \Rhom{R}{{}_{R}M_{S}\Ltensor{S}{}_{S}N}{{}_{R}M_{S}\Ltensor{S}{}_{S}N'}$$ is an isomorphism.

$(3)$ The functor $${}_{R}M_{S}\Ltensor{S}-:\D(S)\rightarrow \D(R)$$ is a full embedding of derived categories.
\eprop

\pf
Condition $(2)$ is just a reformulation of what it means for the functor ${}_{R}M_{S}\Ltensor{S}-$ to be fully faithful.

The functor ${}_{R}M_{S}\Ltensor{S}-:\D(S)\rightarrow \D(R)$ is left adjoint to $$\Rhom{R}{{}_{R}M_{S}}{-}:\D(R)\rightarrow \D(S)$$. The unit of this adjunction is $${}_{S}N\rightarrow \textnormal{RHom}_{R}({}_{R}M_{S},{}_{R}M_{S}\Ltensor{S}{}_{S}N).$$ By \cite[Theorem IV.3.1]{MacLane}, the unit of an adjunction is an isomorphism if and only if the left adjoint is a full embedding. This gives the equivalence $(1)\Leftrightarrow (3)$.
\epf

\prop
\label{new2}
Let $R$ and $S$ be DGAs and suppose that $M$ is a DG $R$-$S$-bimodule. Let $Z=\Rhom{S^{\rm{op}}}{{}_{R}M_{S}}{{}_{S}S_{S}}$ so that $Z$ obtains the structure  ${}_{S}Z_{R}$. Then the following conditions are equivalent:

$(1)$ The canonical map $${}_{S}Z_{R}\Ltensor{R}{}_{R}M_{S}\rightarrow {}_{S}S_{S}$$ is an isomorphism.

$(2)$ For all DG $S$-modules $N$ the canonical map $${}_{S}Z_{R}\Ltensor{R}({}_{R}M_{S}\Ltensor{S}{}_{S}N)\rightarrow {}_{S}N$$ is an isomorphism.

$(3)$ For all DG $S^{\rm{op}}$-modules $N$ and DG $S$-modules $N'$ the canonical map $$(N_{S}\Ltensor{S}{}_{S}Z_{R})\Ltensor{R}({}_{R}M_{S}\Ltensor{S}{}_{S}N')\rightarrow N_{S}\Ltensor{S}{}_{S}N'$$ is an isomorphism.
\eprop

\pf
$(1)\Rightarrow (2).$ Suppose the canonical map ${}_{S}Z_{R}\Ltensor{R}{}_{R}M_{S}\rightarrow {}_{S}S_{S}$ is an isomorphism. Then, for any DG $S$-module $N$ we obtain:
\begin{eqnarray*}
{}_{S}Z_{R}\Ltensor{R}({}_{R}M_{S}\Ltensor{S}{}_{S}N) & \rightiso & ({}_{S}Z_{R}\Ltensor{R}{}_{R}M_{S})\Ltensor{S}{}_{S}N \\
													 & \rightiso & {}_{S}S_{S}\Ltensor{S}{}_{S}N \\													  
 													  &	\rightiso & {}_{S}N,
\end{eqnarray*}
where the first isomorphism is by associativity of left tensor and the second is by hypothesis.

$(2)\Rightarrow (3).$ Suppose that the canonical map ${}_{S}Z_{R}\Ltensor{R}({}_{R}M_{S}\Ltensor{S}{}_{S}N')\rightarrow {}_{S}N'$ is an isomorphism for all DG $S$-modules $N'$. Then we get the following sequence of isomorphisms for DG $S^{\rm{op}}$-modules $N$ and DG $S$-modules $N'$:
\begin{eqnarray*}
(N_{S}\Ltensor{S}{}_{S}Z_{R})\Ltensor{R}({}_{R}M_{S}\Ltensor{S}{}_{S}N') & \rightiso & N_{S}\Ltensor{S}({}_{S}Z_{R}\Ltensor{R}({}_{R}M_{S}\Ltensor{S}{}_{S}N')) \\
	& \rightiso & N_{S}\Ltensor{S}{}_{S}N'.
\end{eqnarray*}

$(3)\Rightarrow (1)$. Suppose that the canonical map 
$$(N_{S}\Ltensor{S}{}_{S}Z_{R})\Ltensor{R}({}_{R}M_{S}\Ltensor{S}{}_{S}N')\rightarrow N_{S}\Ltensor{S}{}_{S}N'$$
 is an isomorphism for all DG $S^{\rm{op}}$-modules $N$ and DG $S$-modules $N'$. Setting $N=N'=S$ we obtain:
\begin{eqnarray*}
{}_{S}Z_{R}\Ltensor{R}{}_{R}M_{S} & \rightiso & ({}_{S}S_{S}\Ltensor{S}{}_{S}Z_{R})\Ltensor{R}({}_{R}M_{S}\Ltensor{S}{}_{S}S_{S}) \\
								  & \rightiso & {}_{S}S_{S}\Ltensor{S}{}_{S}S_{S} \\
								  & \rightiso & {}_{S}S_{S}.
\end{eqnarray*}
This completes the proof.
\epf

\rem
{We now note that if $M$ is finitely built from $S_{S}$ in $\D(S^{\rm{op}})$ then we have the following isomorphisms (c.f. \cite{Recollement}, Setup 2.1):}
\begin{eqnarray*}
\Rhom{S}{{}_{S}Z_{R}}{-} & \cong & \Rhom{S}{\Rhom{S^{\rm{op}}}{{}_{R}M_{S}}{{}_{S}S_{S}}}{-} \\
			 & \cong & {}_{R}M_{S}\Ltensor{S}\Rhom{S}{{}_{S}S_{S}}{-} \\
			 & \cong & {}_{R}M_{S}\Ltensor{S}-.
\end{eqnarray*}
\label{remark1}
\erem

Using Remark \ref{remark1}, we can glue together the conditions of Propositions \ref{new1} and \ref{new2} to link the functor ${}_{R}M_{S}\Ltensor{S} -$ being fully faithful with the conditions of Proposition \ref{new2} which regard the tensor product. This gives us the following generalisation of Theorem 4.4 in \cite{Geigle}.

\thm
\label{new3}
Let $R$ and $S$ be DGAs and suppose $M$ is a DG $R$-$S$-bimodule which is finitely built from $S_{S}$ in $\D(S^{\rm{op}})$.  Let $Z=\Rhom{S^{\rm{op}}}{{}_{R}M_{S}}{{}_{S}S_{S}}$ so that $Z$ obtains the structure ${}_{S}Z_{R}$. Then the following conditions are equivalent:

$(1)$ The canonical map $${}_{S}Z_{R}\Ltensor{R}{}_{R}M_{S}\rightarrow {}_{S}S_{S}$$ is an isomorphism.

$(2)$ For all DG $S$-modules $N$ the canonical map $${}_{S}Z_{R}\Ltensor{R}({}_{R}M_{S}\Ltensor{S}{}_{S}N)\rightarrow {}_{S}N$$ is an isomorphism.

$(3)$ For all DG $S^{\rm{op}}$-modules $N$ and DG $S$-modules $N'$ the canonical map $$(N_{S}\Ltensor{S}{}_{S}Z_{R})\Ltensor{R}({}_{R}M_{S}\Ltensor{S}{}_{S}N')\rightarrow N_{S}\Ltensor{S}{}_{S}N'$$ is an isomorphism.

$(4)$ For all DG $S$-modules $N$ the canonical map $${}_{S}N\rightarrow \Rhom{R}{{}_{R}M_{S}}{{}_{R}M_{S}\Ltensor{S}{}_{S}N}$$ is an isomorphism.

$(5)$ For all DG $S$-modules $N$ and $N'$ the canonical map $$\Rhom{S}{{}_{S}N}{{}_{S}N'}\rightarrow \Rhom{R}{{}_{R}M_{S}\Ltensor{S}{}_{S}N}{{}_{R}M_{S}\Ltensor{S}{}_{S}N'}$$ is an isomorphism.

$(6)$ The functor $${}_{R}M_{S}\Ltensor{S}-:\D(S)\rightarrow \D(R)$$ is a full embedding of derived categories.
\ethm

\pf
The first three conditions are equivalent by Proposition \ref{new2} and conditions $(4)$, $(5)$ and $(6)$ are equivalent by Proposition \ref{new1}. We show that $(1)\Rightarrow (4)$ and $(5)\Rightarrow (2)$.

$(1)\Rightarrow (4).$ Suppose the canonical map ${}_{S}Z_{R}\Ltensor{R}{}_{R}M_{S}\rightarrow {}_{S}S_{S}$ is an isomorphism. Then for all DG $S$-modules $N$ we have the following sequence of isomorphisms.
\begin{eqnarray*}
{}_{S}N & \rightiso & \Rhom{S}{{}_{S}S_{S}}{{}_{S}N} \\
	& \rightiso & \Rhom{S}{{}_{S}Z_{R}\Ltensor{R}{}_{R}M_{S}}{{}_{S}N} \\
	& \rightiso & \Rhom{R}{{}_{R}M_{S}}{\Rhom{S}{{}_{S}Z_{R}}{{}_{S}N}} \\
	& \rightiso & \Rhom{R}{{}_{R}M_{S}}{{}_{R}M_{S}\Ltensor{S}{}_{S}N},
\end{eqnarray*}
where the last isomorphism is given by the duality noted in Remark \ref{remark1}. It can be checked that this composition of maps gives the canonical map.

$(5)\Rightarrow (2)$. It is easy to check that the following diagram commutes for all DG $S$-modules $N$ and $N'$:
$$\xymatrix{\Rhom{S}{N}{N'}\ar[r]^-{\kappa^{*}}\ar[d]_-{\sim}^-{(5)} &   \Rhom{S}{Z\Ltensor{R}(M\Ltensor{S}N)}{N'}\ar[d]^-{\sim}_-{\textnormal{adjunction}}  \\
\Rhom{R}{M\Ltensor{S}N}{M\Ltensor{S}N'}\ar[r]^-{\sim}  & \Rhom{R}{M\Ltensor{S}N}{\Rhom{S}{Z}{N'}}}$$
where $\kappa^{*}$ is the map induced by the canonical map 
$$\kappa:{}_{S}Z_{R}\Ltensor{R}({}_{R}M_{S}\Ltensor{S}{}_{S}N)\rightarrow {}_{S}N,$$
and the bottom arrow is induced by the duality of Remark \ref{remark1}.
Thus $\kappa^{*}$ is an isomorphism for all DG $S$-modules $N$ and $N'$. Fixing $N$, we see that  $$H_{0}(\kappa^{*})=\textnormal{Hom}_{\D(S)}(\kappa, N'):\textnormal{Hom}_{\D(S)}(N,N')\rightarrow \textnormal{Hom}_{\D(S)}(Z\Ltensor{R}(M\Ltensor{S}N),N')$$ is an isomorphism for all DG $S$-modules $N'$. Hence $\kappa$ is an isomorphism for all DG $S$-modules $N$, as required.
\epf

\

Under the additional hypothesis that the DG $R$-$S$-bimodule $M$ is finitely built from ${}_{R}R$ in $\D(R)$ then we have an additional condition equivalent to those in Theorem \ref{new3}.

\prop
\label{new4}
Let $R$ and $S$ be DGAs and $M$ be a DG $R$-$S$-bimodule which is finitely built from ${}_{R}R$ in $\D(R)$. Then the following conditions are equivalent:

\noindent$(1)$ The canonical map ${}_{S}S_{S}\rightarrow \Rhom{R}{{}_{R}M_{S}}{{}_{R}M_{S}}$ is an isomorphism.

\noindent$(2)$ The functor ${}_{R}M_{S}\Ltensor{S}-:\D(S)\rightarrow \D(R)$ is a full embedding of derived categories.
\eprop

\pf
$(2)\Rightarrow (1)$. Suppose that the functor ${}_{R}M_{S}\Ltensor{S}-:\D(S)\rightarrow \D(R)$ is a full embedding of derived categories. By Proposition \ref{new1} the canonical map ${}_{S}N\rightarrow \Rhom{R}{{}_{R}M_{S}}{{}_{R}M_{S}\Ltensor{S}{}_{S}N}$ is an isomorphism for all DG $S$-modules $N$. Setting $N=S$,  we obtain $${}_{S}S_{S}\rightiso  \Rhom{R}{{}_{R}M_{S}}{{}_{R}M_{S}\Ltensor{S}{}_{S}S_{S}}\rightiso \Rhom{R}{{}_{R}M_{S}}{{}_{R}M_{S}}.$$

$(1)\Rightarrow (2)$. Suppose that the canonical map ${}_{S}S_{S}\rightarrow \Rhom{R}{{}_{R}M_{S}}{{}_{R}M_{S}}$ is an isomorphism. Then for any DG $S$-module $N$ we have the following sequence of isomorphisms:
\begin{eqnarray*}
											 {}_{S}N & \rightiso & {}_{S}S_{S}\Ltensor{S}{}_{S}N \\
													 & \rightiso & \Rhom{R}{{}_{R}M_{S}}{{}_{R}M_{S}}\Ltensor{S}{}_{S}N \\
													 & \rightiso & \Rhom{R}{{}_{R}M_{S}}{{}_{R}M_{S}\Ltensor{S}{}_{S}N}
\end{eqnarray*}
where the second isomorphism is by assumption and the third isomorphism is by the fact that $M$ is finitely built from $R$ in $\D(R)$ (c.f. \cite{Recollement}). Thus by Proposition \ref{new1} the functor ${}_{R}M_{S}\Ltensor{S}-:\D(S)\rightarrow \D(R)$ is a full embedding of derived categories.
\epf

\

We can also obtain opposite statements of the above results for DG $S^{\rm{op}}$-modules.

\section{Examples}\label{examples}

We shall now consider three examples of the results in Section \ref{new}. The first example is an application to Dwyer and Greenlees's Morita theory, the second example is an application to Geigle and Lenzing's homological epimorphisms of rings, and the third example is a generalistion of this to homological epimorphisms of DGAs. 

\subsection{Dwyer and Greenlees's Morita theory}\label{dwyer}

The following is taken from J\o rgensen \cite{Recollement} which in turn is taken from Dwyer and Greenlees \cite{Greenlees} with the trivial alteration of replacing ring by DGA.

\setup
{Let $R$ be a DGA and suppose $M$ is a K-projective DG $R$-module which is finitely built from $R$ in $\D(R)$. Let $\F=\textnormal{Hom}_{R}({}_{R}M,{}_{R}M)$ be the endomorphism DGA of $M$, and note that $M$ acquires a left $\F$-structure compatible with its left $R$-structure, ${}_{R,\F}M$. We noted earlier that a DG right  $\F$-module can be canonically identified with a DG left $\Fop$-module (see Section \ref{prelim}). Similarly, a DG left $\F$-module can be canonically identified with a DG right $\Fop$-module. Hence, the left $\F$-structure of $M$ can be regarded as a right $\Fop$-structure, giving $M$ a DG $R$-$\Fop$-structure, ${}_{R}M_{\Fop}$. Note also that, by comparing multiplications, ${}_{\Fop}\F_{\Fop}$ and ${}_{\Fop}(\Fop)_{\Fop}$ are isomorphic as DG $\Fop$-$\Fop$-bimodules.}

\textnormal{Consider the functor:
$$j_{!}(-)={}_{R}M_{\Fop}\Ltensor{\Fop}-:\D(\Fop)\rightarrow \D(R).$$
It is one of the main conclusions of \cite[Section 2]{Greenlees} that this functor is fully faithful. The same fact that $j_{!}(-)$ is a full embedding is a main ingredient in the proof of \cite[Theorem 3.3]{Recollement}. Here we see that full fidelity of $j_{!}(-)$ is also an immediate consequence of Proposition \ref{new4}.} 
\label{setup1}
\esetup 

\prop
\label{dwyer1}
In Setup \ref{setup1} the functor $j_{!}(-):\D(\Fop)\rightarrow \D(R)$ is fully faithful.
\eprop

\pf
Since $M$ is a K-projective  DG left $R$-module we have that 
\begin{eqnarray*}
{}_{\Fop}(\Fop)_{\Fop} & \rightiso & {}_{\Fop}\mathcal{F}_{\Fop} \\
		       & \rightiso & \textnormal{Hom}_{R}({}_{R}M_{\Fop},{}_{R}M_{\Fop}) \\
		       & \rightiso & \Rhom{R}{{}_{R}M_{\Fop}}{{}_{R}M_{\Fop}}. 
\end{eqnarray*}
Since $M$ is finitely built from ${}_{R}R$ in $\D(R)$ it follows by setting $S=\Fop$ in Proposition \ref{new4} that the functor $j_{!}(-)={}_{R}M_{\Fop}\Ltensor{\Fop}- : \D(\Fop)\rightarrow \D(R)$ is a full embedding of derived categories.
\epf 

\

One can also deduce Proposition \ref{dwyer1} directly using Theorem \ref{new3}.

\subsection{Homological epimorphisms of rings}\label{rings}						

We noted in Section \ref{new} that Theorem \ref{new3} was a generalisation of Geigle and Lenzing \cite[Theorem 4.4]{Geigle}. Below we describe how to obtain \cite[Theorem 4.4]{Geigle} from Theorem \ref{new3}.

\rem
{Let $R$ be a ring and let $\D(R)$ be the derived category of $R$. Recall that $$H^{i}(\Rhom{R}{X}{Y})\cong \textnormal{Ext}^{i}_{R}(X,Y)$$ and $$H_{i}(X\Ltensor{R} Y)\cong \textnormal{Tor}^{R}_{i}(X,Y)$$ for two complexes $X$ and $Y$ of $R$-modules, and where Ext and Tor denote the hyperext and hypertor described in \cite{Hartshorne}, respectively. Note also that when $X$ and $Y$ are complexes of $R$-modules concentrated in degree zero the $i^{\textnormal{th}}$-hyperext (resp. $i^{\textnormal{th}}$-hypertor) coincides with the usual Ext (resp. Tor).}
\label{extandtor}
\erem

\rem
{A ring $R$ can be considered as a DGA concentrated in degree zero, which we will also denote by $R$. Moreover, a DG module over $R$ as a DGA is just a complex of $R$-modules, thus the derived categories of $R$ viewed as a ring and $R$ viewed as a DGA coincide.}
\label{ringisdga}
\erem

\setup
{Let $R$ and $S$ be rings and suppose that $\phi: R\rightarrow S$ is a homomorphism of rings. Regard $R$ and $S$ as DGAs concentrated in degree zero; $\phi$ then becomes a morphism of DGAs. Each DG $S$-module can be viewed as a DG $R$-module via $\phi$. In particular, set ${}_{R}M_{S}={}_{R}S_{S}$. Note that ${}_{R}S_{S}$ is a DG $R$-$S$-bimodule which is finitely built from $S_{S}$ in $\D(S^{\rm{o}})$. In the notation of Theorem \ref{new3}, $${}_{S}Z_{R}= \Rhom{S^{\rm{op}}}{{}_{R}S_{S}}{{}_{S}S_{S}}\cong {}_{S}S_{R},$$ so, for instance, in this setup, condition $(1)$ of Theorem \ref{new3} now says that the canonical map  $S\Ltensor{R}S\rightarrow S$ is an isomorphism.}
\label{setup2}
\esetup

We would like to apply Theorem \ref{new3} to the situation of Setup \ref{setup2}. Below we give an indication of how to translate the statements of the conditions in Theorem \ref{new3} to those of \cite[Theorem 4.4]{Geigle} by showing the translation for condition $(1)$.

\newtheorem{misc}[dfn]{Translation}
\begin{misc}
In Setup \ref{setup2}, the canonical map $S\Ltensor{R}S\rightarrow S$ is an isomorphism if and only if the multiplication map $S\otimes_{R}S\rightarrow S$ is an isomorphism and $\textnormal{Tor}_{i}^{R}(S,S)=0$ for all $i\geqslant 1$.
\end{misc}

\pf
Suppose that the canonical map 
$S\Ltensor{R}S\rightarrow S$
is an isomorphism. By Remark \ref{extandtor} we have,  
$$H_{i}(S\Ltensor{R}S)\cong \textnormal{Tor}_{i}^{R}(S,S).$$
Since $S$ is concentrated in degree zero, we have, 
$$\textnormal{Tor}_{i}^{R}(S,S)\cong H_{i}(S\Ltensor{R}S)\cong H_{i}(S) = 0$$
for $i\neq 0$. For $i=0$ we have 
$$S\otimes_{R}S\cong H_{0}(S\Ltensor{R}S)\cong H_{0}(S)=S.$$
It can be seen by an argument involving canonical maps that the composition above is the multiplication map; see Section \ref{canonical}.

Conversely, suppose the multiplication map 
$S\otimes_{R}S\rightarrow S$
is an isomorphism and $\textnormal{Tor}_{i}^{R}(S,S)=0$ for $i\geqslant 1$. 
Since $S$ is concentrated in degree zero, 
$S\Ltensor{R}S$
satisfies 
$H_{i}(S\Ltensor{R}S) = 0$ 
for all $i\leqslant -1$.
For $i\geqslant 1$ we have: 
$$H_{i}(S\Ltensor{R}S)\cong \textnormal{Tor}_{i}^{R}(S,S)=0=H_{i}(S).$$
For $i=0$, we obtain 
$$H_{0}(S\Ltensor{R}S)\cong S\otimes_{R}S\cong S\cong H_{0}(S).$$
Thus 
$S\Ltensor{R}S\rightarrow S$
is a quasi-isomorphism. 
and 
$S$
are objects in a derived category, thus 
$S\Ltensor{R}S\rightarrow S$
is an isomorphism. 
One can again check that this is the canonical map.
\epf

\

Continuing in this manner, we obtain the conditions of the following theorem from the conditions of Theorem \ref{new3}.

\thm
[\cite{Geigle}, Theorem 4.4]
\label{lenzing} 
For a homomorphism of rings $\phi:R\rightarrow S$ the following conditions are equivalent:

$(1)$ The multiplication map ${}_{S}S_{R}\otimes_{R}{}_{R}S_{S}\rightarrow {}_{S}S_{S}$ is an isomorphism and $\textnormal{Tor}^{R}_{i}(S,S)=0$ for all $i\geqslant 1$.

$(2)$ For all left $S$-modules $M$ the multiplication map ${}_{S}S_{R}\otimes_{R}{}_{R}M\rightarrow {}_{S}M$ is an isomorphism and $\textnormal{Tor}^{R}_{i}(S,M)=0$ for all $i\geqslant 1$.

%$(2')$ For all right $S$-modules $N$ the multiplication map $N_{R}\otimes_{R}{}_{R}S_{S}\rightarrow N_{S}$ is an isomorphism and $\textnormal{Tor}^{R}_{i}(N,S)=0$ for all $i\geqslant 1$.

$(3)$ For all right $S$-modules $M$ and all left $S$-modules $N$ the natural map $\textnormal{Tor}^{R}_{i}(M,N)\rightarrow \textnormal{Tor}^{S}_{i}(M,N)$ is an isomorphism for all $i\geqslant 0$.

$(4)$ For all left $S$-modules $M$ the natural map ${}_{S}M\rightarrow \textnormal{Hom}_{R}({}_{R}S_{S},{}_{R}M)$ is an isomorphism and $\textnormal{Ext}^{i}_{R}(S,M)=0$ for all $i\geqslant 1$.

%$(4')$ For all right $S$-modules $N$ the natural map $N_{S}\rightarrow \textnormal{Hom}_{R}({}_{S}S_{R},N_{R})$ is an isomorphism and $\textnormal{Ext}^{i}_{R}(S,N)=0$ for all $i\geqslant 1$.

$(5)$ For all left $S$-modules $M$ and $M'$ the natural map $\textnormal{Ext}^{i}_{S}({}_{S}M,{}_{S}M')\rightarrow\textnormal{Ext}^{i}_{R}({}_{R}M,{}_{R}M')$ is an isomorphism for all $i\geqslant 0$.

%$(5')$ For all right $S$-modules $N$ and $N'$ the natural map $\textnormal{Ext}^{i}_{S}(N_{S},N'_{S})\rightarrow\textnormal{Ext}^{i}_{R}(N_{R},N'_{R})$ is an isomorphism for all $i\geqslant 0$.

$(6)$ The induced functor $\D(\phi):\D(S)\rightarrow \D(R)$ is a full embedding of derived categories.

%$(6')$ The induced functor $\D(\phi^{o}):\D(S^{\rm{op}})\rightarrow \D(R^{o})$ is a full embedding of derived categories.
\ethm

\rem
{It is clear that the corresponding theorem holds for right $S$-modules, and the symmetric nature of condition $(1)$ shows that the statements for right $S$-modules are equivalent to those for left $S$-modules.}
\erem

Theorem \ref{lenzing} as stated above is a strengthening of \cite[Theorem 4.4]{Geigle}. Geigle and Lenzing's theorem is stated in terms of the bounded derived category, whereas the abstract machinery used in Theorem \ref{new3} allows us to work in the full derived category.

\subsection{Homological epimorphisms of DGAs}\label{DGAs}

As in Section \ref{rings}, we obtain the following result from  Theorem \ref{new3} by setting $M=S$. Again, the opposite statements for DG $S^{\rm{op}}$-modules are equivalent to the statements for DG $S$-modules because of the left-right symmetry in statement $(1)$.

\thm
\label{main}
Let $R$ and $S$ be DGAs and let $\phi:R\rightarrow S$ be a morphism of DGAs. 
The following conditions are equivalent:

$(1)$ The canonical map ${}_{S}S_{R}\Ltensor{R}{}_{R}S_{S} \rightarrow {}_{S}S_{S}$ is an isomorphism.

$(2)$ For all DG $S$-modules $M$ the canonical map ${}_{S}S_{R}\Ltensor{R}{}_{R}M\rightarrow {}_{S}M$ is an isomorphism.

%$(2')$ For all DG $S^{\rm{op}}$-modules $N$ the canonical map $N_{R}\Ltensor{R}{}_{R}S_{S}\rightarrow N_{R}$ is an isomorphism.

$(3)$ For all DG $S^{\rm{op}}$-modules $M$ and all DG $S$-modules $N$ the canonical map $M_{R}\Ltensor{R}{}_{R}N\rightarrow M_{S}\Ltensor{S}{}_{S}N$ is an isomorphism.

$(4)$ For all DG $S$-modules $M$ the canonical map ${}_{S}M\rightarrow \Rhom{R}{{}_{R}S_{S}}{{}_{R}M}$ is an isomorphism.

%$(4')$ For all DG $S^{\rm{op}}$-modules $N$ the canonical map $N_{S}\rightarrow \Rhom{R^{o}}{{}_{S}S_{R}}{N_{R}}$ is an isomorphism. 

$(5)$ For all DG $S$-modules $M, M'$ the canonical map $\Rhom{S}{{}_{S}M}{{}_{S}M'}\rightarrow \Rhom{R}{{}_{R}M}{{}_{R}M'}$ is an isomorphism.

%$(5')$ For all DG $S^{\rm{op}}$-modules $N, N'$ the canonical map $\Rhom{S^{\rm{op}}}{N_{S}}{N'_{S}}\rightarrow \Rhom{R^{o}}{N_{R}}{N'_{R}}$ is an isomorphism.

$(6)$ The induced functor $\D(\phi):\D(S)\rightarrow\D(R)$ is a full embedding of derived categories.

%$(6')$ The induced functor $\D(\phi^{o}):\D(S^{\rm{op}})\rightarrow\D(R^{o})$ is a full embedding of derived categories.
\ethm

Following Geigle and Lenzing, \cite{Geigle}, we make the following definition.
 
\df
{We call a morphism of DGAs $\phi:R\rightarrow S$ a \textit{homological epimorphism of DGAs} if it satisfies the equivalent conditions of Theorem \ref{main}.}
\edf

For a morphism $\phi: R\rightarrow S$ of DGAs, if $S$ is finitely built from $R$ as a DG $R$-module, then we can add an extra condition equivalent to the conditions in Theorem \ref{main}.

\prop
\label{prop1}
Let $\phi:R\rightarrow S$ be a morphism of DGAs. Suppose that $S$ is finitely built from $R$ as a DG $R$-module. Then the following conditions are equivalent:

$(1)$ The canonical map ${}_{S}S_{S}\rightarrow \Rhom{R}{{}_{R}S_{S}}{{}_{R}S_{S}}$ is an isomorphism.

$(2)$ The morphism $\phi:R\rightarrow S$ is a homological epimorphism.
\eprop

\pf
Follows from Proposition \ref{new4}.
\epf

\

In Theorem \ref{main} a homological epimorphism of DGAs is characterised by the canonical map 
$S\Ltensor{R}S\rightarrow S$
being an isomorphism. However, the strengthened condition that $S$ is finitely built from ${}_{R}R$ in $\D(R)$ is required to prove that a homological epimorphism of DGAs can also be characterised by the canonical map 
$S\rightarrow \Rhom{R}{S}{S}$
being an isomorphism. This apparent lack of symmetry arises because the characterisation of a homological epimorphism of DGAs involving left tensor requires the trivial fact that $S$ is finitely built from ${}_{S}S$ in $\D(S)$; whereas, the characterisation with RHom requires that $S$ is finitely built from ${}_{R}R$ in $\D(R)$ (c.f. Theorem \ref{new3} and Proposition \ref{new4}).

\

\noindent
\textbf{Acknowledgement.}  I would like to thank my supervisor, Peter J\o rgensen, for all the help and advice he has given during the preparation of this paper, and also to thank the University of Leeds and EPSRC of the United Kingdom for financial support.

\end{document}